# ANOTHER SET OF SEQUENCES, SUB-SEQUENCES, AND SEQUENCES OF SEQUENCES

by Florentin Smarandache, Ph. D.
University of New Mexico
Gallup, NM 87301, USA

**Abstract.** New sequences in number theory are showed below with definitions, examples, solved or open questions and references for each case.

**Introduction.**
In this paper 101 new integer sequences, sub-sequences, and sequences of sequences, together with related unsolved problems and conjectures, are presented.

**Sequences of sequences**:

1. THE DIGIT SEQUENCES.

General definition:
in any numeration base B, for any given infinite integer or rational sequence $S_1, S_2, S_3, \ldots,$ and any digit D from 0 to B-1,
it's built up a new integer sequence witch
associates to $S_1$ the number of digits D of $S_1$ in base B,
to $S_2$ the number of digits D of $S_2$ in base B, and so on...

For exemple, considering the prime number sequence in base 10, then the number of digits 1 (for exemple) of each prime number following their order is: 0 0 0 0 2 1 1 1 0 0 1 0...
(The digit-1 prime sequence)

Second exemple if we consider the factorial sequence n! in base 10, then the number of digits 0 of each factorial number following their order is: 0 0 0 0 0 1 1 2 2 1 3...
(The digit-0 factorial sequence)

Third exemple if we consider the sequence n^n in base 10, n=1,2,..., then the number of digits 5 of each term 1^1, 2^2, 3^3,..., following their order is: 0 0 0 1 1 1 1 0 0 0...
(The digit-5 n^n sequence)

## 2. THE CONSTRUCTION SEQUENCE.

General definition:

in any numeration base B, for any given infinite integer or rational sequence $S_1$, $S_2$, $S_3$, ..., and any digits $D_1$, $D_2$, ..., $D_k$ ($k < B$),
it's built up a new integer sequence such that
each of its terms $Q_1 < Q_2 < Q_3 < \ldots$ is formed by these digits $D_1$, $D_2$, ..., $D_k$ only (all these digits are used), and matches a term $S_i$ of the previous sequence.

For exemple, considering in base 10 the prime number sequence,
And, say, digits 1 and 7,
we construct a written-only-with-these-digits (all these digits are used) prime number new sequence: 17 71 ...
(The digit-1-7-only prime sequence)

Second exemple, considering in base 10 the multiple of 3 sequence,
and the digits 0 and 1,
we construct a written-only-with-these-digits (all these digits are used) multiple of 3 new sequence: 1011 1101 1110 10011 10101 10110 11001 11010 11100 ...
(The digit-0-1-only multiple of 3 sequence)

3. THE CONSECUTIVE SEQUENCE:

1 12 123 1234 12345 123456 1234567 12345678 123456789 12345678910 1234567891011 123456789101112 12345678910111213 ...

How many primes are there among these numbers?

In a general form, the Consecutive Sequence is considered in an arbitrary numeration base B.

4. THE SYMMETRIC SEQUENCE:

1 11 121 1221 12321 123321 1234321 12344321 123454321 1234554321 12345654321 123456654321 1234567654321 12345677654321 123456787654321 1234567887654321 12345678987654321 123456789987654321 12345678910987654321 1234567891010987654321 12345678910110987654321 1234567891011110987654321 123456789101112110987654321 12345678910111212110987654321 1234567891011121312110987654321 ...

How many primes are there among these numbers?

In a general form, the Symmetric Sequence is considered in an arbitrary numeration base B.

5)THE GENERAL RESIDUAL SEQUENCE:

$(x + C_1)\ldots(x + C_{F(m)})$, $m = 2, 3, 4, \ldots$,

where $C_i$, $1 \le i \le F(m)$, forms a reduced set of residues mod $m$, $x$ is an integer, and $F$ is Euler's totient.

The Smarandache General Residual Sequence is induced from the The Smarandache Residual Function (see <Libertas Mathematica>):

Let $L : Z \times Z \to Z$ be a function defined by

$L(x,m)=(x + C_1)\ldots(x + C_{F(m)})$,

$_1$ $_{F(m)}$

where $C_i$, $1 <= i <= F(m)$, forms a reduced set of residues mod m, m >= 2, x is an integer, and F is Euler's totient.

The Smarandache Residual Function is important because it generalixes the classical theorems by Wilson, Fermat, Euler, Wilson, Gauss, Lagrange, Leibnitz, Moser, and Sierpinski all together.

For x=0 it's obtained the following sequence:

$L(m) = C_1 \dots C_{F(m)}$, where m = 2, 3, 4, ...

(the product of all residues of a reduced set mod m):

1 2 3 24 5 720 105 2240 189 3628800 385 479001600 19305 896896 2027025 20922789888000 85085 6402373705728000 8729721 47297536000 1249937325 ...

which is found in "The Enciclopedia of Integer Sequences".

The Residual Function extends it.

## 6)THE NUMERICAL CARPET:

has the general form

```
              .
              .
              .
              1
             1a1
            1aba1
           1abcba1
          1abcdcba1
         1abcdedcba1
        1abcdefedcba1
   ...1abcdefgfedcba1...
        1abcdefedcba1
         1abcdedcba1
          1abcdcba1
           1abcba1
            1aba1
             1a1
              1
              .
              .
              .
```

```
On the border of level 0, the elements are equal to "1";
  they form a rhomb.
Next, on the border of level 1, the elements are equal to "a",
  where "a" is the sum of all elements of the previous border;
  the "a"s form a rhomb too inside the previous one.
Next again, on the border of level 2, the elements are equal to "b",
  where "b" is the sum of all elements of the previous border;
  the "b"s form a rhomb too inside the previous one.
And so on...
The carpet is symmetric and esthetic, in its middle g is the
sum of all carpet numbers (the core).

Look at a few terms of the Numerical Carpet:

                              1

                              1
                             141
                              1

                              1
                           1  8  1
                        1  8 40  8  1
                           1  8  1
                              1

                              1
                          1  12   1
                      1  12 108  12  1
                   1 12 108 540 108 12 1
                      1  12 108  12  1
                          1  12   1
                              1

                              1
                         1   16    1
                    1   16  208   16   1
                1  16  208 1872  208  16  1
             1 16 208 1872 9360 1872 208 16 1
                 1 16  208 1872  208  16  1
                    1   16  208   16   1
                         1   16    1
                              1

                              1
                       1     20     1
                 1    20    340    20    1
            1   20   340   4420   340   20   1
        1  20  340  4420  39780  4420  340  20  1
     1 20 340 4420 39780 198900 39780 4420 340 20 1
        1  20  340  4420  39780  4420  340  20  1
            1   20   340   4420   340   20   1
                 1    20    340    20    1
                       1     20     1
                              1
```

.
.
.

Or, under another form:

| | | | | | | | | | |
|---|---|---|---|---|---|---|---|---|---|
| 1 | | | | | | | | | |
| 1 | 4 | | | | | | | | |
| 1 | 8 | 40 | | | | | | | |
| 1 | 12 | 108 | 504 | | | | | | |
| 1 | 16 | 208 | 1872 | 9360 | | | | | |
| 1 | 20 | 340 | 4420 | 39780 | 198900 | | | | |
| 1 | 24 | 504 | 8568 | 111384 | 1002456 | 5012280 | | | |
| 1 | 28 | 700 | 14700 | 249900 | 3248700 | 29238300 | 146191500 | | |
| 1 | 32 | 928 | 23200 | 487200 | 8282400 | 107671200 | 969040800 | 4845204000 | |

..............................................................
.
.
.

General Formula:

$$C(n,k) = 4n\prod_{i=1}^{k}(4n-4i+1) \text{ for } 1 \le k \le n,$$

and $C(n,0) = 1$.

7)THE SQUARE COMPLEMENTS:

1 2 3 1 5 6 7 2 1 10 11 3 14 15 1 17 2 19 5 21 22 23 6 1 26 3 7 29 30 31 2 33 34 35 1 37 38 39 10 41 42 43 11 5 46 47 3 1 2 51 13 53 6 55 14 57 58 59 15 61 62 7 1 65 66 67 17 69 70 71 2 ...

Definition:

for each integer n to find the smallest integer k such that nk is a perfect square..

(All these numbers are square free.)

8)THE CUBIC COMPLEMENTS:

1 4 9 2 25 36 49 1 3 100 121 18 169 196 225 4 289 12 361 50 441 484 529
9 5 676 1 841 900 961 2 1089 1156 1225 6 1369 1444 1521 25 1681 1764 1849
242 75 2116 2209 36 7 20 ...

Definition:
for each integer n to find the smallest integer k such that
nk is a perfect cub.
(All these numbers are cub free.)

9)THE M-POWER COMPLEMENTS (generalization):

Definition:
for each integer n to find the smallest integer k such that
nk is a perfect m-power (m => 2).
(All these numbers are m-power free.)

10)THE SQUARE FREE SIEVE:

2 3 5 6 7 10 11 13 14 15 17 19 21 22 23 26 29 30 31 33 34 35 37 38 39 41
42 43 46 47 51 53 55 57 58 59 61 62 65 66 67 69 70 71 ...

Definition: from the set of natural numbers (except 0 and 1):
- take off all multiples of 2^2 (i.e. 4, 8, 12, 16, 20, ...)
- take off all multiples of 3^2
- take off all multiples of 5^2

... and so on (take off all multiples of all square primes).

(One obtains all square free numbers.)

11)THE CUBE FREE SIEVE:

2 3 4 5 6 7 9 10 11 12 13 14 15 17 18 19 20 21 22 23 25 26 28 29 30 31 33
34 35 36 37 38 39 41 42 43 44 45 46 47 49 50 51 52 53 55 57 58 59 60 61 62
63 65 66 67 68 69 70 71 73 ...

Definition: from the set of natural numbers (except 0 and 1):
- take off all multiples of 2^3 (i.e. 8, 16, 24, 32, 40, ...)
- take off all multiples of 3^3
- take off all multiples of 5^3

... and so on (take off all multiples of all cubic primes).

(One obtains all cube free numbers)

12)THE M-POWER FREE SIEVE (generalization):

Definition: from the set of natural numbers (except 0 and 1) take off all multiples of 2^m, afterwards all multiples of 3^m, ... and so on (take off all multiples of all m-power primes, m >= 2).

(One obtains all m-power free numbers.)

14)THE IRRATIONAL ROOT SIEVE:

2 3 5 6 7 10 11 12 13 14 15 17 18 19 20 21 22 23 24 26 28 29 30 31 33 34 35 37 38 39 40 41 42 43 44 45 46 47 48 50 51 52 53 54 55 56 57 58 59 60 61 62 63 65 66 67 68 69 70 71 72 73 ...

Definition: from the set of natural numbers (except 0 and 1):

- take off all powers of 2^k, k >= 2, (i.e. 4, 8, 16, 32, 64, ...)
- take off all powers of 3^k, k >= 2;
- take off all powers of 5^k, k >= 2;
- take off all powers of 6^k, k >= 2;
- take off all powers of 7^k, k >= 2;
- take off all powers of 10^k, k >= 2;

... and so on (take off all k-powers, k >= 2, of all square free numbers -- see the square free sieve).

(One obtains all natural numbers those m-th roots, for any m >= 2, are irrational.)

14)THE SYLLABIC PUZZLE:

1 1 1 1 1 1 2 2 2 1 3 1 3 3 3 3 4 3 4 2 ...

(a(n) = the number of syllables of n in English language).

15)THE CODE PUZZLE:

151405 202315 2008180505 06152118 06092205 190924 1905220514 0509070820 14091405 200514 051205220514 ...

Using the following letter-to-number code:

A B C D E F G H I J K L M N O P Q R S T U V W X Y Z

01 02 03 04 05 06 07 08 09 10 11 12 13 14 15 16 17 18 19 20 21 22 23 24 25 26
then a(n) = the numerical code for the spelling of n in English
language;  for exemple: 1 = ONE = 151405, etc.

16)THE PIERCED CHAIN:

```
101 1010101 10101010101 101010101010101 1010101010101010101
10101010101010101010101 101010101010101010101010101 ...
(a(n) = 101 * 1 0001 0001 ... 0001 , for n >= 1)
                |  | |  | ... |  |
                ---- ----     ----
                  1    2       n-1
```

How many a(n)/101 are primes ?

17) THE Smarandache QUOTIENTS:

1 1 2 6 24 1 720 3 80 12 3628800 2 479001600 360 8 45 20922789888000
40 6402373705728000 6 240 1814400 1124000727777607680000 1 145152
239500800 13440 180 304888344611713860501504000000 ...
(For each n to find the smallest k such that nk is a factorial number.)

18)THE (INFERIOR) PRIME PART:

2 3 3 5 5 7 7 7 7 11 11 13 13 13 13 17 17 19 19 19 19 23 23 23 23 23 23
29 29 31 31 31 31 31 31 37 37 37 37 41 41 43 43 43 43 47 47 47 47 47 47
53 53 53 53 53 53 59 ...
(For any positive real number n one defines a(n) as the largest prime number less than or equal to n)

19)THE (SUPERIOR) PRIME PART:

2 2 2 3 5 5 7 7 11 11 11 11 13 13 17 17 17 17 19 19 23 23 23 23 29 29 29
29 29 29 31 31 37 37 37 37 37 37 41 41 41 41 43 43 47 47 47 47 53 53 53
53 53 53 59 59 59 59 59 59 61 ...

(For any positive real number n one defines a(n) as the smallest prime number greater than or equal to n)

20)THE DOUBLE FACTORIAL COMPLEMENTS:

1 1 1 2 3 8 15 1 105 192 945 4 10395 46080 1 3 2027025 2560 34459425 192
5 3715891200 13749310575 2 81081 1961990553600 35 23040 213458046676875
128 6190283353629375 12 ...

(For each n to find the smallest k such that nk is a double factorial,
i.e. nk = either 1*3*5*7*9*...*n if n is odd,
either 2*4*6*8*...*n if n is even.)

21)THE PRIME COMPLEMENTS:

1 0 0 1 0 1 0 3 2 1 0 1 0 3 2 1 0 1 0 3 2 1 0 5 4 3 2 1 0 1 0 5 4 3 2 1 0
3 2 1 0 1 0 3 2 1 0 5 4 3 2 1 0 ...

(For each n to find the smallest k such that n+k is prime.)

Remark: Is it possible to get as large as we want but finite decreasing sequence k, k-1, k-2, ..., 2, 1, 0 (odd k) included in the previous sequence -- i.e. for any even integer are there two primes those difference is equal to it? I conjecture the answer is negative.

22)THE ROMANIAN LETTERS ORDER:

E A I R T S P O C U D Z N L V M F G B H X J K W Q Y

(The Romanian language letters frequency in the juridical texts according to a study done by F. Smarandache)

23)THE ODD SIEVE:

7 13 19 23 25 31 33 37 43 47 49 53 55 61 63 67 73 75 79 83 85 91 93 97 ...

(All odd numbers that are not equal to the difference of two primes.)

A sieve is used to get this sequence:

- substract 2 from all prime numbers and obtain a temporary sequence;
- choose all odd numbers that do not belong to the temporary one.

24)THE DOUBLE FACTORIAL NUMBERS:

1 2 3 4 5 6 7 4 9 10 11 6 13 14 5 6 17 12 19 10 7 22 23 6 15 26 9 14 29 10 31 8 11 34 7 12 37 38 13 10 41 14 43 22 9 46 47 6 21 10 ...

(a(n) is the smallest integer such that a(n)!! is a multiple of n.)

25)THE Smarandache PARADOXIST NUMBERS:

There exist a few "Smarandache" number sequences.
A number n is called a "Smarandache paradoxist number" if and only if n doesn't belong to any of the S defined numbers.

Question: find the Smarandcahe paradoxist number sequence.

Solution?

If a number k is a S paradoxist number, then k doesn't belong to any of the Smarandache defined numbers,
therefore k doesn't belong to the Smarandache paradoxist numbers too!

If a number k doesn't belong to any of the Smarandache defined numbers, then k is a Smarandache paradoxist number,
therefore k belongs to a Smarandache defined numbers (because Smarandache paradoxist numbers is also in the same category) -- contradiction.

Dilemma: is the Smarandache paradoxist number sequence empty ??

26)THE NON-SMARANDACHE NUMBERS:

A number n is called a "non-Smarandache number" if and only if n is neither a Smarandache paradoxist number nor any of the Smarandache defined numbers.

Question: find the non-Smarandache number sequence.

Dilemma 1: is the non-Smarandache number sequence empty, too ??

Dilemma 2: is a non-Smarandache number equivalent to a Smarandache paradoxist number ??? (this would be another paradox !! ... because a non-Smarandache number is not a Smarandache paradoxist number).

27)THE PARADOX OF SMARANDACHE NUMBERS:

Any number is a Smarandache number, the non-Smarandache number too.
(This is deduced from the following paradox (see the reference):
"All is possible, the impossible too!")

28) Prime base:

0 1 10 100 101 1000 1001 10000 10001 10010 10100 100000 100001 1000000 1000001 1000010 1000100 10000000 10000001 100000000 100000001 100000010 100000100 1000000000 1000000001 1000000010 1000000100 1000000101 ...
(Each number n written into the prime base.)

(One defines over the set of natural numbers the following infinite base: $p_0 = 1$, and for $k \geq 1$ $p_k$ is the k-th prime number.)

He proved that every positive integer A may be uniquely written into the prime base as:

$$A = (\overline{a \ldots a\ a}) \stackrel{def}{===} \sum^{n} a\ p$$ , with all a = 0 or 1, (of course a = 1),

```
      n      1 0 (SP)      /   i i            i                     n
                           ---
                           i=0
```

in the following way:

- if $p_n \le A < p_{n+1}$ then $A = p_n + r_1$;
- if $p_m \le r_1 < p_{m+1}$ then $r_1 = p_m + r_2$, $m < n$;

and so on until one obtains a rest $r_j = 0$.

Therefore, any number may be written as a sum of prime numbers + e, where e = 0 or 1.

If we note by p(A) the superior part of A (i.e. the largest prime less than or equal to A), then A is written into the prime base as:

$$A = p(A) + p(A-p(A)) + p(A-p(A)-p(A-p(A))) + ...$$

This base is important for partitions with primes.

29) Deconstructive sequence:

```
1 23 456 7891 23456 789123 4567891 23456789 123456789 1234567891  ...
|          ||         ||        ||        | |       | |       ||
 ----------  ---------  --------  --------   --------  -------  - ...
```

30)Goldbach-Smarandache table:

6 10 14 18 26 30 38 42 42 54 62 74 74 90 ...

(a(n) is the largest even number such that any other even number not exceeding it is the sum of two of the first n odd primes.)

It helps to better understand Goldbach's conjecture:

- if a(n) is unlimited, then the conjecture is true;
- if a(n) is constant after a certain rank, then the conjecture is false.

Also, the table gives how many times an even number is written as a sum of two odd primes, and in what combinations -- which can be found in the "Encyclopedia of Integer Sequences" by N. J. A. Sloane and S. Plouffe, Academic Press, 1995.

Of course, $a(n) <= 2p_n$ , where $p_n$ is the n-th odd prime, n = 1, 2, 3, ... .

Here is the table:

| + | 3 | 5 | 7 | 11 | 13 | 17 | 19 | 23 | 29 | 31 | 37 | 41 | 43 | 47 | . . . |
|---|---|---|---|---|---|---|---|---|---|---|---|---|---|---|---|
| 3 | 6 | 8 | 10 | 14 | 16 | 20 | 22 | 26 | 32 | 34 | 40 | 44 | 46 | 50 | . |
| 5 | | 10 | 12 | 16 | 18 | 22 | 24 | 28 | 34 | 36 | 42 | 46 | 48 | 52 | . |
| 7 | | | 14 | 18 | 20 | 24 | 26 | 30 | 36 | 38 | 44 | 48 | 50 | 54 | . |
| 11 | | | | 22 | 24 | 28 | 30 | 34 | 40 | 42 | 48 | 52 | 54 | 58 | . |
| 13 | | | | | 26 | 30 | 32 | 36 | 42 | 44 | 50 | 54 | 56 | 60 | . |
| 17 | | | | | | 34 | 36 | 40 | 46 | 48 | 54 | 58 | 60 | 64 | . |
| 19 | | | | | | | 38 | 42 | 48 | 50 | 56 | 60 | 62 | 66 | . |
| 23 | | | | | | | | 46 | 52 | 54 | 60 | 64 | 66 | 70 | . |
| 29 | | | | | | | | | 58 | 60 | 66 | 70 | 72 | 76 | . |
| 31 | | | | | | | | | | 62 | 68 | 72 | 74 | 78 | . |
| 37 | | | | | | | | | | | 74 | 78 | 80 | 84 | . |
| 41 | | | | | | | | | | | | 82 | 84 | 88 | . |
| 43 | | | | | | | | | | | | | 86 | 90 | . |
| 47 | | | | | | | | | | | | | | 94 | . |
| . . . | . | . | . | . | . | . | . | . | . | . | . | . | . | . | . |

31) Primitive numbers (of power 2):
2 4 4 6 8 8 8 10 12 12 14 16 16 16 16 18 20 20 22 24 24 24 26 28 28 30 32 32 32 32 32 34 36 36 38 40 40 40 42 44 44 46 48 48 48 48 50 52 52 54 56 56 56 58 60 60 62 64 64 64 64 64 64 66 ...
(a(n) is the smallest integer such that a(n)! is divisible by 2^n)

Curious property: this is the sequence of even numbers, each number being repeated as many times as its exponent (of power 2) is.

This is one of irreductible functions, noted $S_2(k)$, which helps to calculate the Smarandache function (called also Smarandache numbers in "The Encyclopedia of Integer Sequences", by N. J. A. Sloane and S. Plouffe, Academic Press, 1995).

32) Primitive numbers (of power 3):
3 6 9 9 12 15 18 18 21 24 27 27 27 30 33 36 36 39 42 45 45 48 51 54 54 54 57 60 63 63 66 69 72 72 75 78 81 81 81 81 84 87 90 90 93 96 99 99 102 105 108 108 108 111 ...
(a(n) is the smallest integer such that a(n)! is divisible by 3^n)

Curious property: this is the sequence of multiples of 3, each number being repeated as many times as its exponent (of power 3) is.

This is one of irreductible functions, noted $S_3(k)$, which helps to calculate the function (called also numbers in "The Encyclopedia of Integer Sequences", by N. J. A. Sloane and S. Plouffe, Academic Press, 1995).

33)Primitive numbers (of power p, p prime) -- generalization:
(a(n) is the smallest integer such that a(n)! is divisible by p^n)

Curious property: this is the sequence of multiples of p, each number being repeated as many times as its exponent (of power p) is.

These are the irreductible functions, noted $S_p(k)$, for any prime number p, which helps to calculate the function (called also numbers in "The Encyclopedia of Integer Sequences", by N. J. A. Sloane and S. Plouffe, Academic Press, 1995).

34) Square residues:
1 2 3 2 5 6 7 2 3 10 11 6 13 14 15 2 17 6 19 10 21 22 23 6 5 26 3 14 29 30 31 2 33 34 35 6 37 38 39 10 41 42 43 22 15 46 47 6 7 10 51 26 53 6 14 57 58 59 30 61 62 21 ...
(a(n) is the largest square free number which divides n.)

Or, a(n) is the number n released of its squares:
if $n = (p_1 \wedge a_1) * ... * (p_r \wedge a_r)$, with all $p_i$ primes and all $a_i >= 1$,
then $a(n) = p_1 * ... * p_r$.

Remark: at least the (2^2)*k-th numbers (k = 1, 2, 3, ...) are released of their squares;
and more general: all (p^2)*k-th numbers (for all p prime, and k = 1, 2, 3, ...) are released of their squares.

35) Cubical residues:
1 2 3 4 5 6 7 4 9 10 11 12 13 14 15 4 17 18 19 20 21 22 23 12 25 26 9 28 29 30 31 4 33 34 35 36 37 38 39 20 41 42 43 44 45 46 47 12 49 50 51 52 53 18 55 28 ...
(a(n) is the largest cube free number which divides n.)

Or, a(n) is the number n released of its cubicals:
if $n = (p_1 \wedge a_1) * ... * (p_r \wedge a_r)$, with all $p_i$ primes and all $a_i >= 1$,
then $a(n) = (p_1 \wedge b_1) * ... * (p_r \wedge b_r)$, with all $b_i = \min \{2, a_i\}$.

Remark: at least the (2^3)*k-th numbers (k = 1, 2, 3, ...) are released of their cubicals;
and more general: all (p^3)*k-th numers (for all p prime, and k = 1, 2, 3, ...) are released of their cubicals.

36) m-power residues (generalization):
(a(n) is the largest m-power free number which divides n.)

Or, a(n) is the number n released of its m-powers:
if $n = (p_1 ^{a_1}) * \ldots * (p_r ^{a_r})$, with all $p_i$ primes and all $a_i >= 1$,
then $a(n) = (p_1 ^{b_1}) * \ldots * (p_r ^{b_r})$, with all $b_i = \min \{ m-1, a_i \}$.

Remark: at least the (2^m)*k-th numbers (k = 1, 2, 3, ...) are released of their m-powers;
and more general: all (p^m)*k-th numers (for all p prime, and k = 1, 2, 3, ...) are released of their m-powers.

37) Exponents (of power 2):
0 1 0 2 0 1 0 3 0 1 0 2 0 1 0 4 0 1 0 2 0 1 0 2 0 1 0 2 0 1 0 5 0 1 0 2 0
1 0 3 0 1 0 2 0 1 0 3 0 1 0 2 0 1 0 3 0 1 0 2 0 1 0 6 0 1 ...
(a(n) is the largest exponent (of power 2) which divides n)

Or, a(n) = k if 2^k divides n but 2^(k+1) does not.

38) Exponents (of power 3):
0 0 1 0 0 1 0 0 2 0 0 1 0 0 1 0 0 2 0 0 1 0 0 1 0 0 3 0 0 1 0 0 1 0 0 2 0
0 1 0 0 1 0 0 2 0 0 1 0 0 1 0 0 2 0 0 1 0 0 1 0 0 2 0 0 1 0 ...
(a(n) is the largest exponent (of power 3) which divides n)

Or, a(n) = k if 3^k divides n but 3^(k+1) does not.

39) Exponents (of power p) -- generalization :
(a(n) is the largest exponent (of power p) which divides n,
where p is an integer >= 2)

Or, a(n) = k if p^k divides n but p^(k+1) does not.

(reviewed in <Zentralblatt fur Mathematik> by P. Kiss: 11002, pre744, 1992;
and <The American Mathematical Monthly>, Aug.-Sept. 1991);
Arizona State University, Hayden Library, "The Florentin Smarandache papers" special collection, Tempe, AZ 85287-1006, USA, phone: (602)965-6515 (Carol Moore librarian), email: ICCLM@ASUACAD.BITNET .

40) Pseudo-primes of first kind:
2,3,5,7,11,13,14,16,17,19,20,23,29,30,31,32,34,35,37,38,41,43,47,50,53,59,
61,67,70,71,73,74,76,79,83,89,91,92,95,97,98,101,103,104,106,107,109,110,
112,113,115,118,119,121,124,125,127,128,130,131,133,134,136,137,139,140,
142,143,145,146, ...
(A number is a pseudo-prime of first kind if some permutation of the digits is a prime number, including the identity permutation.)

(Of course, all primes are pseudo-primes of first kind, but not the reverse!)

41)Pseudo-primes of second kind:
14,16,20,30,32,34,35,38,50,70,74,76,91,92,95,98,104,106,110,112,115,118,
119,121,124,125,128,130,133,134,136,140,142,143,145,146, ...
(A composite number is a pseudo-prime of second kind if some permutation of the digits is a prime number.)

42)Pseudo-primes of third kind:
11,13,14,16,17,20,30,31,32,34,35,37,38,50,70,71,73,74,76,79,91,92,95,97,98,
101,103,104,106,107,109,110,112,113,115,118,119,121,124,125,127,128,130,
131,133,134,136,137,139,140,142,143,145,146, ...
(A number is a pseudo-prime of third kind if some nontrivial permutation of the digits is a prime number.)

Question: How many pseudo-primes of third kind are prime numbers? (he conjectured: an infinity).
(There are primes which are not pseudo-primes of third kind, and the reverse:
there are pseudo-primes of third kind which are not primes.)

43)Ppseudo-squares of first kind:
1,4,9,10,16,18,25,36,40,46,49,52,61,63,64,81,90,94,100,106,108,112,121,136,
144,148,160,163,169,180,184,196,205,211,225,234,243,250,252,256,259,265,
279,289,295,297,298,306,316,324,342,360,361,400,406,409,414,418,423,432,
441,448,460,478,481,484,487,490,502,520,522,526,529,562,567,576,592,601,
603,604,610,613,619,625,630,631,640,652,657,667,675,676,691,729,748,756,
765,766,784,792,801,810,814,829,841,844,847,874,892,900,904,916,925,927,
928,940,952,961,972,982,1000, ...
(A number is a pseudo-square of first kind if some permutation of the digits is a perfect square, including the identity permutation.)

(Of course, all perfect squares are pseudo-squares of first kind, but not the reverse!)

One listed all pseudo-squares of first kind up to 1000.

44) Pseudo-squares of second kind:
10,18,40,46,52,61,63,90,94,106,108,112,136,148,160,163,180,184,205,211,234,
243,250,252,259,265,279,295,297,298,306,316,342,360,406,409,414,418,423,
432,448,460,478,481,487,490,502,520,522,526,562,567,592,601,603,604,610,
613,619,630,631,640,652,657,667,675,691,748,756,765,766,792,801,810,814,
829,844,847,874,892,904,916,925,927,928,940,952,972,982,1000, ...
(A non-square number is a pseudo-square of second kind if some permutation of the digits is a square.)

One listed all pseudo-squares of second kind up to 1000.

45) Pseudo-squares of third kind:
10,18,40,46,52,61,63,90,94,100,106,108,112,121,136,144,148,160,163,169,180,
184,196,205,211,225,234,243,250,252,256,259,265,279,295,297,298,306,316,
342,360,400,406,409,414,418,423,432,441,448,460,478,481,484,487,490,502,520,
522,526,562,567,592,601,603,604,610,613,619,625,630,631,640,652,657,667,
675,676,691,748,756,765,766,792,801,810,814,829,844,847,874,892,900,904,
916,925,927,928,940,952,961,972,982,1000,...
(A number is a pseudo-square of third kind if some nontrivial permutation of the digits is a square.)

Question: How many pseudo-squares of third kind are square numbers? (he conjectured: an infinity).
(There are squares which are not pseudo-squares of third kind, and the reverse:
there are pseudo-squares of third kind which are not squares.)

One listed all pseudo-squares of third kind up to 1000.

46) Pseudo-cubes of first kind:
1,8,10,27,46,64,72,80,100,125,126,152,162,207,215,216,251,261,270,279,297,
334,343,406,433,460,512,521,604,612,621,640,702,720,729,792,800,927,972,
1000,...
(A number is a pseudo-cube of first kind if some permutation of the digits is a cube, including the identity permutation.)

(Of course, all perfect cubes are pseudo-cubes of first kind, but not the reverse!)

One listed all pseudo-cubes of first kind up to 1000.

47) Pseudo-cubes of second kind:
10,46,72,80,100,126,152,162,207,215,251,261,270,279,297,334,406,433,460,
521,604,612,621,640,702,720,792,800,927,972,...
(A non-cube number is a pseudo-cube of second kind if some permutation of the digits is a cube.)

One listed all pseudo-cubes of second kind up to 1000.

48) Pseudo-cubes of third kind:
10,46,72,80,100,125,126,152,162,207,215,251,261,270,279,297,334,343,
406,433,460,512,521,604,612,621,640,702,720,792,800,927,972,1000,...
(A number is a pseudo-cube of third kind if some nontrivial permutation of the digits is a cube.)

Question: How many pseudo-cubes of third kind are cubes? (he conjectured: an infinity).
(There are cubes which are not pseudo-cubes of third kind, and the reverse: there are pseudo-cubes of third kind which are not cubes.)

One listed all pseudo-cubes of third kind up to 1000.

49)Pseudo-m-powers of first kind:
(A number is a pseudo-m-power of first kind if some permutation of the digits is an m-power, including the identity permutation; m >= 2.)

50)Pseudo-m-powers of second kind:
(An m-power number is a pseudo-m-power of second kind if some permutation of the digits is an m-power; m >= 2.)

51) Pseudo-m-powers of third kind:
(A number is a pseudo-m-power of third kind if some nontrivial permutation of the digits is an m-power; m >= 2.)

Question: How many pseudo-m-powers of third kind are m-power numbers? (he conjectured: an infinity).
(There are m-powers which are not pseudo-m-powers of third kind, and the reverse: there are pseudo-m-powers of third kind which are not m-powers.)

S. Plouffe, Academic Press, 1995;
also online, email: superseeker@research.att.com ( SUPERSEEKER by N. J. A. Sloane, S. Plouffe, B. Salvy, ATT Bell Labs, Murray Hill, NJ 07974, USA);

52) Mirror sequence:

1 212 32123 4321234 543212345 65432123456 7654321234567 876543212345678 98765432123456789 109876543212345678910 1110987654321234567891011 ...

Question: How many of them are primes?

53) Permutation sequence:

12 1342 135642 13578642 13579108642 135791112108642 1357911131412108642 13579111315161412108642 135791113151718161412108642 1357911131517192018161412108642 ...

Question: Is there any perfect power among these numbers?
(Their last digit should be:
either 2 for exponents of the form 4k+1,
either 8 for exponents of the form 4k+3, where k >= 0 .)

I conjecture: no!

54)Generalizated permutation sequence:

If g(n), as a function, gives the number of digits of a(n), and F if a permutation of g(n) elements, then:

$a(n) = \overline{F(1)F(2)...F(g(n))}$ .

55)Constructive set (of digits 1,2):

1,2,11,12,21,22,111,112,121,122,211,212,221,222,1111,1112,1121,1122,1211, 1212,1221,1222,21112112,2121,2122,2211,2212,2221,2222,...
(Numbers formed by digits 1 and 2 only.)

Definition:

a1) 1, 2 belong to S;

a2) if a, b belong to S, then $\overline{ab}$ belongs to S too;

a3) only elements obtained by rules a1) and a2) applied a finite number of times belong to S.

Remark:

- there are 2^k numbers of k digits in the sequence, for k = 1, 2, 3, ... ;
- to obtain from the k-digits number group the (k+1)-digits number group, just put first the digit 1 and second the digit 2 in the front of all k-digits numbers.

56)Constructive set (of digits 1,2,3):

1,2,3,11,12,13,21,22,23,31,32,33,111,112,113,121,122,123,131,132,133,211, 212,213,221,222,223,231,232,233,311,312,313,321,322,323,331,332,333,...

(Numbers formed by digits 1, 2, and 3 only.)

Definition:

a1) 1, 2, 3 belong to S;

a2) if a, b belong to S, then $\overline{ab}$ belongs to S too;

a3) only elements obtained by rules a1) and a2) applied a finite number of times belong to S.

Remark:

- there are 3^k numbers of k digits in the sequence, for k = 1, 2, 3, ... ;
- to obtain from the k-digits number group the (k+1)-digits number group, just put first the digit 1, second the digit 2, and third the digit 3 in the front of all k-digits numbers.

57)Generalizated constructive set:

(Numbers formed by digits $d_1$, $d_2$, ..., $d_m$ only, all $d_i$ being different each other, 1 <= m <= 9.)

Definition:

a1) $d_1$, $d_2$, ..., $d_m$ belong to S;

a2) if a, b belong to S, then $\overline{ab}$ belongs to S too;

a3) only elements obtained by rules a1) and a2) applied a finite number of times belong to S.

Remark:

- there are m^k numbers of k digits in the sequence, for k = 1, 2, 3, ... ;
- to obtain from the k-digits number group the (k+1)-digits number group, just put first the digit $d_1$, second the digit $d_2$, ..., and the m-th time digit $d_m$ in the front of all k-digits numbers.

More general: all digits $d_i$ can be replaced by numbers as large as we want (therefore of many digits each), and also m can be as large as we want.

58)Square roots:

0,1,1,1,2,2,2,2,2,3,3,3,3,3,3,3,4,4,4,4,4,4,4,4,4,5,5,5,5,5,5,5,5,5,5, 6,6,6,6,6,6,6,6,6,6,6,6,6,7,7,7,7,7,7,7,7,7,7,7,7,7,7,7,8,8,8,8,8,8,8,8, 8,8,8,8,8,8,8,8,9,9,9,9,9,9,9,9,9,9,9,9,9,9,9,9,9,9,9,10,10,10,10,10,10,10,

10,10,10,10,10,10,10,10,10,10,10,10,10,10,...
(a(n) is the superior integer part of square root of n.)

Remark: this sequence is the natural sequence, where each number is repeated 2n+1 times,
because between n^2 (included) and (n+1)^2 (excluded) there are (n+1)^2 - n^2 different numbers.

59)Cubical roots:
0,1,1,1,1,1,1,1,2,2,2,2,2,2,2,2,2,2,2,2,2,2,2,2,2,2,2,3,3,3,3,3,3,3,3,3,3,
3,3,3,3,3,3,3,3,3,3,3,3,3,3,3,3,3,3,3,3,3,3,3,3,3,3,3,4,4,4,4,4,4,4,4,4,4,
4,4,4,4,4,4,4,4,4,4,4,4,4,4,4,4,4,4,4,4,4,4,4,4,4,4,4,4,4,4,4,4,4,4,4,4,4,
4,4,4,4,4,4,4,4,4,4,4,4,4,4,...
(a(n) is the superior integer part of cubical root of n.)

Remark: this sequence is the natural sequence, where each number is repeated 3n^2 + 3n + 1 times,
because between n^3 (included) and (n+1)^3 (excluded) there are (n+1)^3 - n^3 different numbers.

60) m-power roots:
(a(n) is the superior integer part of m-power root of n.)

Remark: this sequence is the natural sequence, where each number is repeated (n+1)^m - n^m times.

61)Pseudo-factorials of first kind:
1,2,6,10,20,24,42,60,100,102,120,200,201,204,207,210,240,270,402,420,600,
702,720,1000,1002,1020,1200,2000,2001,2004,2007,2010,2040,2070,2100,2400,
2700,4002,4005,4020,4050,4200,4500,5004,5040,5400,6000,7002,7020,7200,...
(A number is a pseudo-factorial of first kind if some permutation of the digits is a factorial number, including the identity permutation.)

(Of course, all factorials are pseudo-factorials of first kind, but not the reverse!)

One listed all pseudo-factorials of first kind up to 10000.

Procedure to obtain this sequence:
- calculate all factorials with one digit only (1!=1, 2!=2, and 3!=6), this is line_1 (of one digit pseudo-factorials):
  1,2,6;
- add 0 (zero) at the end of each element of line_1, calculate all factorials with two digits (4!=24 only) and all permutations of their digits:
  this is line_2 (of two digits pseudo-factorials):
  10,20,60; 24, 42;
- add 0 (zero) at the end of each element of line_2 as well as anywhere in between their digits,

calculate all factorials with three digits (5!=120, and 6!=720)
and all permutations of their digits:
this is line_3 (of three digits pseudo-factorials):
100,200,600,240,420,204,402; 120,720, 102,210,201,702,270,720;

and so on ...
to get from line_k to line_(k+1) do:
- add 0 (zero) at the end of each element of line_k as well as anywhere in between their digits,
calculate all factorials with (k+1) digits
and all permutations of their digits;

The set will be formed by all line_1 to the last line elements in an increasing order.

The pseudo-factorials of second kind and third kind can be deduced from the first kind ones..

62)Pseudo-factorials of second kind:

10,20,42,60,100,102,200,201,204,207,210,240,270,402,420,600,
702,1000,1002,1020,1200,2000,2001,2004,2007,2010,2040,2070,2100,2400,
2700,4002,4005,4020,4050,4200,4500,5004,5400,6000,7002,7020,7200,...

(A non-factorial number is a pseudo-factorial of second kind if some permutation of the digits is a factorial number.)

63)Pseudo-factorials of third kind:

10,20,42,60,100,102,200,201,204,207,210,240,270,402,420,600,
702,1000,1002,1020,1200,2000,2001,2004,2007,2010,2040,2070,2100,2400,
2700,4002,4005,4020,4050,4200,4500,5004,5400,6000,7002,7020,7200,...

(A number is a pseudo-factorial of third kind if some nontrivial permutation of the digits is a factorial number.)

Question: How many pseudo-factorials of third kind are factorial numbers? (he conjectured: none! ... that means the pseudo-factorials of second kind set and pseudo-factorials of third kind set coincide!).

64)Pseudo-divisors of first kind:

1,10,100,1,2,10,20,100,200,1,3,10,30,100,300,1,2,4,10,20,40,100,200,400,
1,5,10,50,100,500,1,2,3,6,10,20,30,60,100,200,300,600,1,7,10,70,100,700,
1,2,4,8,10,20,40,80,100,200,400,800,1,3,9,10,30,90,100,300,900,1,2,5,10,
20,50,100,200,500,1000,...

(The pseudo-divisors of first kind of n)

(A number is a pseudo-divisor of first kind of n if some permutation of the digits is a divisor of n, including the identity permutation.)

(Of course, all divisors are pseudo-divisors of first kind, but not the reverse!)

A strange property: any integer has an infinity of pseudo-divisors of first kind !!
because 10...0 becomes 0...01 = 1, by a circular permutation of its digits, and 1 divides any integer !

One listed all pseudo-divisors of first kind up to 1000 for the numbers 1, 2, 3, ..., 10.

Procedure to obtain this sequence:

- calculate all divisors with one digit only, this is line_1 (of one digit pseudo-divisors);
- add 0 (zero) at the end of each element of line_1, calculate all divisors with two digits and all permutations of their digits: this is line_2 (of two digits pseudo-divisors);
- add 0 (zero) at the end of each element of line_2 as well as anywhere in between their digits, calculate all divisors with three digits and all permutations of their digits: this is line_3 (of three digits pseudo-divisors);

and so on ...

to get from line_k to line_(k+1) do:

- add 0 (zero) at the end of each element of line_k as well as anywhere in between their digits, calculate all divisors with (k+1) digits and all permutations of their digits;

The set will be formed by all line_1 to the last line elements in an increasing order.

The pseudo-divisors of second kind and third kind can be deduced from the first kind ones.

65) Pseudo-divisors of second kind:

10,100,10,20,100,200,10,30,100,300,10,20,40,100,200,400,10,50,100,500,10,
20,30,60,100,200,300,600,10,70,100,700,10,20,40,80,100,200,400,800,10,30,
90,100,300,900,20,50,100,200,500,1000,...

(The pseudo-divisors of second kind of n)

(A non-divisor of n is a pseudo-divisor of second kind of n if some permutation of the digits is a divisor of n.)

66)Pseudo-divisors of third kind:

10,100,10,20,100,200,10,30,100,300,10,20,40,100,200,400,10,50,100,500,10,
20,30,60,100,200,300,600,10,70,100,700,10,20,40,80,100,200,400,800,10,30,
90,100,300,900,10,20,50,100,200,500,1000,...

(The pseudo-divisors of third kind of n)

(A number is a pseudo-divisor of third kind of n if some nontrivial permutation of the digits is a divisor of n.)

A strange property: any integer has an infinity of pseudo-divisors of third kind !! because 10...0 becomes 0...01 = 1, by a circular permutation of its digits, and 1 divides any integer !

There are divisors of n which are not pseudo-divisors of third kind of n, and the reverse: there are pseudo-divisors of third kind of n which are not divisors of n.

67)Pseudo-even numbers of first kind:
0,2,4,6,8,10,12,14,16,18,20,21,22,23,24,25,26,27,28,29,30,32,34,36,38,40,
41,42,43,44,45,46,47,48,49,50,52,54,56,58,60,61,62,63,64,65,66,67,68,69,70,
72,74,76,78,80,81,82,83,84,85,86,87,88,89,90,92,94,96,98,100,...
(The pseudo-even numbers of first kind)

(A number is a pseudo-even number of first kind if some permutation of the digits is a even number, including the identity permutation.)

(Of course, all even numbers are pseudo-even numbers of first kind, but not the reverse!)

A strange property: an odd number can be a pseudo-even number!

One listed all pseudo-even numbers of first kind up to 100.

68)Pseudo-even numbers of second kind:
21,23,25,27,29,41,43,45,47,49,61,63,65,67,69,81,83,85,87,89,101,103,105,
107,109,121,123,125,127,129,141,143,145,147,149,161,163,165,167,169,181,
183,185,187,189,201,...
(The pseudo-even numbers of second kind)

(A non-even number is a pseudo-even number of second kind if some permutation of the digits is a even number.)

69)Pseudo-even numbers of third kind:
20,21,22,23,24,25,26,27,28,29,40,41,42,43,44,45,46,47,48,49,60,61,62,63,64,
65,66,67,68,69,80,81,82,83,84,85,86,87,88,89,100,101,102,103.104,105,106,
107,108,109,110,120,121,122,123,124,125,126,127,128,129,130,...
(The pseudo-even numbers of third kind)

(A number is a pseudo-even number of third kind if some nontrivial permutation of the digits is a even number.)

70)Pseudo-multiples of first kind (of 5):
0,5,10,15,20,25,30,35,40,45,50,51,52,53,54,55,56,57,58,59,60,65,70,75,80,
85,90,95,100,101,102,103,104,105,106,107,108,109,110,115,120,125,130,135,
140,145,150,151,152,153,154,155,156,157,158,159,160,165,...
(The pseudo-multiples of first kind of 5)

(A number is a pseudo-multiple of first kind of 5 if some permutation of the digits is a multiple of 5, including the identity permutation.)

(Of course, all multiples of 5 are pseudo-multiples of first kind, but not the reverse!)

71)Pseudo-multiples of second kind (of 5):
51,52,53,54,56,57,58,59,101,102,103,104,106,107,108,109,151,152,153,154,
156,157,158,159,201,202,203,204,206,207,208,209,251,252,253,254,256,257,
258,259,301,302,303,304,306,307,308,309,351,352...
(The pseudo-multiples of second kind of 5)

(A non-multiple of 5 is a pseudo-multiple of second kind of 5 if some permutation of the digits is a multiple of 5.)

72)Pseudo-multiples of third kind (of 5):
50,51,52,53,54,55,56,57,58,59,100,101,102,103,104,105,106,107,108,109,110,
115,120,125,130,135,140,145,150,151,152,153,154,155,156,157,158,159,160,
165,170,175,180,185,190,195,200,...
(The pseudo-multiples of third kind of 5)

(A number is a pseudo-multiple of third kind of 5 if some nontrivial permutation of the digits is a multiple of 5.)

Generalizations:

73)Pseudo-multiples of first kind of p (p is an integer >= 2):
(The pseudo-multiples of first kind of p)

(A number is a pseudo-multiple of first kind of p if some permutation of the digits is a multiple of p, including the identity permutation.)

(Of course, all multiples of p are pseudo-multiples of first kind, but not the reverse!)

Procedure to obtain this sequence:
- calculate all multiples of p with one digit only (if any), this is line_1 (of one digit pseudo-multiples of p);
- add 0 (zero) at the end of each element of line_1, calculate all multiples of p with two digits (if any) and all permutations of their digits: this is line_2 (of two digits pseudo-multiples of p);
- add 0 (zero) at the end of each element of line_2 as well as anywhere in between their digits, calculate all multiples with three digits (if any) and all permutations of their digits: this is line_3 (of three digits pseudo-multiples of p);

and so on ...
to get from line_k to line_(k+1) do:
- add 0 (zero) at the end of each element of line_k as well as anywhere in between their digits, calculate all multiples with (k+1) digits (if any) and all permutations of their digits;

The set will be formed by all line_1 to the last line elements in an increasing order.

The pseudo-multiples of second kind and third kind of p can be deduced from the first kind ones.

74)Pseudo-multiples of second kind of p (p is an integer >= 2):
(The pseudo-multiples of second kind of p)

(A non-multiple of p is a pseudo-multiple of second kind of p if some permutation of the digits is a multiple of p.)

75)Pseudo-multiples of third kind of p (p is an integer >= 2):
(The pseudo-multiples of third kind of p)

(A number is a pseudo-multiple of third kind of p if some nontrivial permutation of the digits is a multiple of p.)

76)Binary sieve:
1,3,5,9,11,13,17,21,25,27,29,33,35,37,43,49,51,53,57,59,65,67,69,73,75,77, 81,85,89,91,97,101,107,109,113,115,117,121,123,129,131,133,137,139,145, 149,...
(Starting to count on the natural numbers set at any step from 1:
- delete every 2-nd numbers
- delete, from the remaining ones, every 4-th numbers
... and so on: delete, from the remaining ones, every (2^k)-th numbers, k = 1, 2, 3, ... .)

Conjectures:
- there are an infinity of primes that belong to this sequence;
- there are an infinity of numbers of this sequence which are not prime.

77)Trinary sieve:
1,2,4,5,7,8,10,11,14,16,17,19,20,22,23,25,28,29,31,32,34,35,37,38,41,43,46, 47,49,50,52,55,56,58,59,61,62,64,65,68,70,71,73,74,76,77,79,82,83,85,86,88, 91,92,95,97,98,100,101,103,104,106,109,110,112,113,115,116,118,119,122,124, 125,127,128,130,131,133,137,139,142,143,145,146,149,...
(Starting to count on the natural numbers set at any step from 1:
- delete every 3-rd numbers
- delete, from the remaining ones, every 9-th numbers
... and so on: delete, from the remaining ones, every (3^k)-th numbers, k = 1, 2, 3, ... .)

Conjectures:
- there are an infinity of primes that belong to this sequence;
- there are an infinity of numbers of this sequence which are not prime.

78) n-ary sieve (generalization, n >= 2):

(Starting to count on the natural numbers set at any step from 1:
- delete every n-th numbers
- delete, from the remaining ones, every (n^2)-th numbers
- ... and so on: delete, from the remaining ones, every (n^k)-th numbers, k = 1, 2, 3, ... .)

Conjectures:
- there are an infinity of primes that belong to this sequence;
- there are an infinity of numbers of this sequence which are not prime.

79)Consecutive sieve:

1,3,5,9,11,17,21,29,33,41,47,57,59,77,81,101,107,117,131,149.153,173,191,
209,213,239,257,273,281,321,329,359,371,401,417,441,435,491,...

(From the natural numbers set:
- keep the first number,
  delete one number out of 2 from all remaining numbers;
- keep the first remaining number,
  delete one number out of 3 from the next remaining numbers;
- keep the first remaining number,
  delete one number out of 4 from the next remaining numbers;

... and so on, for step k (k >= 2):
- keep the first remaining number,
  delete one number out of k from the next remaining numbers;

... .)

This sequence is much less dense than the prime number sequence, and their ratio tends to $p_n : n$ as n tends to infinity.

For this sequence we chosen to keep the first remaining number at all steps, but in a more general case: the kept number may be any among the remaining k-plet (even at random).

80)General sequence-sieve:

Let $u_i > 1$, for i = 1, 2, 3, ..., a strictly increasing positive integer sequence. Then:

From the natural numbers set:
- keep one number among 1, 2, 3, ..., $u_1 - 1$,
  and delete every $u_1$-th numbers;
- keep one number among the next $u_2 - 1$ remaining numbers,

and delete every $u_2$-th numbers;

... and so on, for step k (k >= 1):

- keep one number among the next $u_k$ - 1 remaining numbers,

  and delete every $u_k$-th numbers;

...

Problem: study the relationship between sequence $u_i$, i = 1, 2, 3, ..., and the remaining sequence resulted from the general sieve.

$u_i$, previously defined, id called sieve generator.

81) (Inferior) square part:

0,1,1,1,4,4,4,4,4,9,9,9,9,9,9,9,16,16,16,16,16,16,16,16,16,25,25,25,25,25, 25,25,25,25,25,25,36,36,36,36,36,36,36,36,36,36,36,36,36,49,49,49,49,49, 49,49,49,49,49,49,49,49,49,49,64,64,...

(The largest square less than or equal to n)

82) (Superior) square part:

0,1,4,4,4,9,9,9,9,9,16,16,16,16,16,16,16,25,25,25,25,25,25,25,25,25,36,36, 36,36,36,36,36,36,36,36,36,49,49,49,49,49,49,49,49,49,49,49,49,49,64,64,64, 64,64,64,64,64,64,64,64,64,64,64,64,81,81,...

(The smallest square greater than or equal to n)

83)(Inferior) cube part:

0,1,1,1,1,1,1,1,8,8,8,8,8,8,8,8,8,8,8,8,8,8,8,8,8,8,8,27,27,27,27,27,27, 27,27,27,27,27,27,27,27,27,27,27,27,27,27,27,27,27,27,27,27,27,27,27,27, 27,27,27,27,27,27,27,64,64,64,...

(The largest cube less than or equal to n)

84)(Superior) cube part:

0,1,8,8,8,8,8,8,8,27,27,27,27,27,27,27,27,27,27,27,27,27,27,27,27,27,27, 27,64,64,64,64,64,64,64,64,64,64,64,64,64,64,64,64,64,64,64,64,64,64,64, 64,64,64,64,64,64,64,64,64,64,64,64,64,64,125,125,125,..

(The smallest cube greater than or equal to n)

85)(Inferior) factorial part:

1,2,2,2,2,6,6,6,6,6,6,6,6,6,6,6,6,6,6,6,6,6,6,24,24,24,24,24,24,24,24,24, 24,24,24,24,24,24,24,24,24,...

(a(n) is the largest factorial less than or equal to n.)

86) (Superior) factorial part:

1,2,6,6,6,6,24,24,24,24,24,24,24,24,24,24,24,24,24,24,24,24,24,24,120,120,

120,120,120,120,120,120,120,120,120,...
(a(n) is the smallest factorial greater than or equal to n.)

87)Digital sum:

```
0,1,2,3,4,5,6,7,8,9,1,2,3,4,5,6,7,8,9,10,2,3,4,5,6,7,8,9,10,11,
|                 | |                  | |                   |
 -----------------   -----------------   -------------------

3,4,5,6,7,8,9,10,11,12,4,5,6,7,8,9,10,11,12,13,5,6,7,8,9,10,11,12,13,14,...
|                    | |                     | |                      |
 -------------------   --------------------   ---------------------
```

(a(n) is the sum of digits.)

88)Digital products:

```
0,1,2,3,4,5,6,7,8,9,0,1,2,3,4,5,6,7,8,9,0,2,4,6,8,19,12,14,16,18,
|                 | |                | |                      |
 -----------------   ----------------   ---------------------

0,3,6,9,12,15,18,21,24,27,0,4,8,12,16,20,24,28,32,36,0,5,10,15,20,25,...
|                       | |                        | |
 -----------------------   -----------------------   -------------- ...
```

(a(n) is the product of digits.)

89)Divisor products:
1,2,3,8,5,36,7,64,27,100,11,1728,13,196,225,1024,17,5832,19,8000,441,484,
23,331776,125,676,729,21952,29,810000,31,32768,1089,1156,1225,10077696,37,
1444,1521,2560000,41,...
(a(n) is the product of all positive divisors of n)

90)Proper divisor products:
1,1,1,2,1,6,1,8,3,10,1,144,1,14,15,64,1,324,1,400,21,22,1,13824,5,26,27,
784,1,27000,1,1024,33,34,35,279936,1,38,39,64000,1,...
(a(n) is the product of all positive divisors of n but n)

91)Pseudo-odd numbers of first kind:
1,3,5,7,9,10,11,12,13,14,15,16,17,18,19,21,23,25,27,29,30,31,32,33,34,35,
36,37,38,39,41,43,45,47,49,50,51,52,53,54,55,56,57,58,59,61,63,65,67,69,70,
71,72,73,74,75,76,...
(Some permutation of digits is an odd number)

92)Pseudo-odd numbers of second kind:
10,12,14,16,18,30,32,34,36,38,50,52,54,56,58,70,72,74,76,78,90,92,94,96,98,
100,102,104,106,108,110,112,114,116,118,...
(Even numbers such that some permutation of digits is an odd number)

93)Pseudo-odd numbers of third kind:

10,11,12,13,14,15,16,17,18,19,30,31,32,33,34,35,36,37,38,39,50,51,52,53,54,55,56,57,58,59,70,71,72,73,74,75,76,...

(Nontrivial permutation of digits is an odd number)

94)Pseudo-triangular numbers:

1,3,6,10,12,15,19,21,28,30,36,45,54,55,60,61,63,66,78,82,87,91,...

(Some permutation of digits is a triangular number)

A triangular number has the general form: $n(n+1)/2$.

95)Square base:

0,1,2,3,10,11,12,13,20,100,101,102,103,110,111,112,1000,1001,1002,1003,1010,1011,1012,1013,1020,10000,10001,10002,10003,10010,10011,10012,10013,10020,10100,10101,100000,100001,100002,100003,100010,100011,100012,100013,100020,100100,100101,100102,100103,100110,100111,100112,101000,101001,101002,101003,101010,101011,101012,101013,101020,101100,101101,101102,1000000,...

(Each number n written into the square base.)

(One defines over the set of natural numbers the following infinite base: for $k \ge 0$ $s_k = k^2$.)

He proved that every positive integer A may be uniquely written into the square base as:

$$A = \overline{(a_n \ldots a_1 a_0)}_{(S2)} \overset{def}{=} \sum_{i=0}^{n} a_i s_i \text{, with } a_i = 0 \text{ or } 1 \text{ for } i \ge 2,$$

$0 \le a_0 \le 3$, $0 \le a_1 \le 2$, and of course $a_n = 1$,
in the following way:

- if $s_n \le A < s_{n+1}$ then $A = s_n + r_1$;
- if $s_m \le r_1 < p_{m+1}$ then $r_1 = s_m + r_2$, $m < n$;

and so on until one obtains a rest $r_j = 0$.

Therefore, any number may be written as a sum of squares (1 not counted as a square -- being obvious) + e, where e = 0, 1, or 3.

If we note by s(A) the superior square part of A (i.e. the largest square less than or equal to A), then A is written into the square base as:

$$A = s(A) + s(A-s(A)) + s(A-s(A)-s(A-s(A))) + \ldots$$

This base is important for partitions with squares.

96) m-power base (generalization):
(Each number n written into the m-power base, where m is an integer >= 2.)

(One defines over the set of natural numbers the following infinite m-power base: for $k \ge 0$ $t_k = k^m$.)

He proved that every positive integer A may be uniquely written into the m-power base as:

$$A = \overline{(a_n \ldots a_1 a_0)}_{(SM)} \overset{def}{=} \sum_{i=0}^{n} a_i t_i, \text{ with } a_i = 0 \text{ or } 1 \text{ for } i \ge m,$$

$$0 \le a_i \le \left[ ((i+2)^m - 1) / (i+1)^m \right] \text{ (integer part)}$$

for $i = 0, 1, \ldots, m-1$, $a_i = 0$ or $1$ for $i \ge m$, and of course $a_n = 1$,

in the following way:

- if $t_n \le A < t_{n+1}$ then $A = t_n + r_1$;
- if $t_m \le r_1 < t_{m+1}$ then $r_1 = t_m + r_2$, $m < n$;

and so on until one obtains a rest $r_j = 0$.

Therefore, any number may be written as a sum of m-powers (1 not counted as an m-power -- being obvious) + e, where e = 0, 1, 2, ..., or 2^m-1.

If we note by t(A) the superior m-power part of A (i.e. the largest m-power less than or equal to A), then A is written into the m-power base as:

$$A = t(A) + t(A-t(A)) + t(A-t(A)-t(A-t(A))) + \ldots$$

This base is important for partitions with m-powers.

97) Generalized base:
(Each number n written into the generalized base.)

(One defines over the set of natural numbers the following infinite generalized base: $1 = g_0 < g_1 < \ldots < g_k < \ldots$ .)

He proved that every positive integer A may be uniquely written into the generalized base as:

$$\overline{A = (a_n \ldots a_1 a_0)}_{(SG)} \overset{def}{=} \sum_{i=0}^{n} a_i g_i \text{, with } 0 \le a_i \le \left[ (g_{i+1} - 1) / g_i \right]$$

(integer part) for i = 0, 1, ..., n, and of course $a_n \ge 1$,

in the following way:

- if $g_n \le A < g_{n+1}$ then $A = g_n + r_1$ ;
- if $g_m \le r_1 < g_{m+1}$ then $r_1 = g_m + r_2$ , $m < n$;

and so on untill one obtains a rest $r_j = 0$.

If we note by g(A) the superior generalized part of A (i.e. the largest $g_i$ less than or equal to A), then A is written into the m-power base as:

$$A = g(A) + g(A-g(A)) + g(A-g(A)-g(A-g(A))) + \ldots$$

This base is important for partitions: the generalized base may be any infinite integer set (primes, squares, cubes, any m-powers, Fibonacci/Lucas numbers, Bernoully numbers, Smarandache numbers, etc.) those partitions are studied.

A particular case is when the base verifies: $2g_i \ge g_{i+1}$ for any i, and $g_0 = 1$, because all coefficients of a written number into this base will be 0 or 1.

98)Smarandache-Vinogradov table:

9,15,21,29,39,47,57,65,71,93,99,115,129,137,...

(a(n) is the largest odd number such that any odd number >= 9 not exceeding it is the sum of three of the first n odd primes.)

It helps to better understand Goldbach's conjecture for three primes:

- if a(n) is unlimited, then the conjecture is true;
- if a(n) is constant after a certain rank, then the conjecture is false.

(Vinogradov proved in 1937 that any odd number greater than 3^(3^15) satisfaies this conjecture.
But what about values less than 3^(3^15) ?)

Also, the table gives you in how many different combinations an odd number is written as a sum of three odd primes, and in what combinations.

Of course, $a(n) \le 3p_n$ , where $p_n$ is the n-th odd prime, n = 1, 2, 3, ... .

It is also generalized for the sum of m primes,
and how many times a number is written as a sum of m primes (m > 2).

This is a 3-dimensional 14x14x14 table, that we can expose only as 14 planar 14x14 tables (using Goldbach-Smarandache table):

| 3 + | 3 | 5 | 7 | 11 | 13 | 17 | 19 | 23 | 29 | 31 | 37 | 41 | 43 | 47 | . . . |
|---|---|---|---|---|---|---|---|---|---|---|---|---|---|---|---|
| 3 | 9 | 11 | 13 | 17 | 19 | 23 | 25 | 29 | 35 | 37 | 43 | 47 | 49 | 53 | . |
| 5 | | 13 | 15 | 19 | 21 | 25 | 27 | 31 | 37 | 39 | 45 | 49 | 51 | 55 | . |
| 7 | | | 17 | 21 | 23 | 27 | 29 | 33 | 39 | 41 | 47 | 51 | 53 | 57 | . |
| 11 | | | | 25 | 27 | 31 | 33 | 37 | 43 | 45 | 51 | 55 | 57 | 61 | . |
| 13 | | | | | 29 | 33 | 35 | 39 | 45 | 47 | 53 | 57 | 59 | 63 | . |
| 17 | | | | | | 37 | 39 | 43 | 49 | 51 | 57 | 61 | 63 | 67 | . |
| 19 | | | | | | | 41 | 45 | 51 | 53 | 59 | 63 | 65 | 69 | . |
| 23 | | | | | | | | 49 | 55 | 57 | 63 | 67 | 69 | 73 | . |
| 29 | | | | | | | | | 61 | 63 | 69 | 73 | 75 | 79 | . |
| 31 | | | | | | | | | | 65 | 71 | 75 | 77 | 81 | . |
| 37 | | | | | | | | | | | 77 | 81 | 83 | 87 | . |
| 41 | | | | | | | | | | | | 85 | 87 | 91 | . |
| 43 | | | | | | | | | | | | | 89 | 93 | . |
| 47 | | | | | | | | | | | | | | 97 | . |
| | . . . | | | | | | | | | | | | | | . . . |

| 5 + | 3 | 5 | 7 | 11 | 13 | 17 | 19 | 23 | 29 | 31 | 37 | 41 | 43 | 47 | . . . |
|---|---|---|---|---|---|---|---|---|---|---|---|---|---|---|---|
| 3 | 11 | 13 | 15 | 19 | 21 | 25 | 27 | 31 | 37 | 39 | 45 | 49 | 51 | 55 | . |
| 5 | | 15 | 17 | 21 | 23 | 27 | 29 | 33 | 39 | 41 | 47 | 51 | 53 | 57 | . |
| 7 | | | 19 | 23 | 25 | 29 | 31 | 35 | 41 | 43 | 49 | 53 | 55 | 59 | . |
| 11 | | | | 27 | 29 | 33 | 35 | 39 | 45 | 47 | 53 | 57 | 59 | 63 | . |
| 13 | | | | | 31 | 35 | 37 | 41 | 47 | 49 | 55 | 59 | 61 | 65 | . |
| 17 | | | | | | 39 | 41 | 45 | 51 | 53 | 59 | 63 | 65 | 69 | . |
| 19 | | | | | | | 43 | 47 | 53 | 55 | 61 | 65 | 67 | 71 | . |
| 23 | | | | | | | | 51 | 57 | 59 | 65 | 69 | 71 | 75 | . |
| 29 | | | | | | | | | 63 | 65 | 71 | 75 | 77 | 81 | . |
| 31 | | | | | | | | | | 67 | 73 | 77 | 79 | 83 | . |
| 37 | | | | | | | | | | | 79 | 83 | 85 | 89 | . |
| 41 | | | | | | | | | | | | 87 | 89 | 93 | . |
| 43 | | | | | | | | | | | | | 91 | 95 | . |
| 47 | | | | | | | | | | | | | | 99 | . |
| | . . . | | | | | | | | | | | | | | . . . |

```
 -----
| 7   |
|  +  |
|     | 3  5  7 11 13 17 19 23 29 31 37 41 43  47
 ----- -------------------------------------------- . . .
    3 |13 15 17 21 23 27 29 33 39 41 47 51 53  57 .
    5 |   17 19 23 25 29 31 35 41 43 49 53 55  59 .
    7 |      21 25 27 31 33 37 43 45 51 55 57  61 .
   11 |         29 31 35 37 41 47 49 55 59 61  65 .
   13 |            33 37 39 43 49 51 57 61 63  67 .
   17 |               41 43 47 53 55 61 65 67  71 .
   19 |                  45 49 55 57 63 67 69  73 .
   23 |                     53 59 61 67 71 73  77 .
   29 |                        65 67 73 77 79  83 .
   31 |                           69 75 79 81  85 .
   37 |                              81 85 87  91 .
   41 |                                 89 91  95 .
   43 |                                    93  97 .
   47 |                                       101 .
      ...............................................
      .                                               .
      .                                                 .
      .                                                   .


 -----
|11   |
|  +  |
|     | 3  5  7 11 13 17 19 23 29 31 37 41 43  47
 ----- -------------------------------------------- . . .
    3 |17 19 21 25 27 31 33 37 43 45 51 55 57  61 .
    5 |   21 23 27 29 33 35 39 45 47 53 57 59  63 .
    7 |      25 29 31 35 37 41 47 49 55 59 61  65 .
   11 |         33 35 39 41 45 51 53 59 63 65  69 .
   13 |            37 41 43 47 53 55 61 65 67  71 .
   17 |               45 47 51 57 59 65 69 71  75 .
   19 |                  49 53 59 61 67 71 73  77 .
   23 |                     57 63 65 71 75 77  81 .
   29 |                        69 71 77 81 83  87 .
   31 |                           73 79 83 85  89 .
   37 |                              85 89 91  95 .
   41 |                                 93 95  99 .
   43 |                                    97 101 .
   47 |                                       105 .
      ...............................................
      .                                               .
      .                                                 .
      .                                                   .
```

```
 -----
|13   |
|  +  |
|     | 3  5  7 11 13 17 19 23 29 31 37 41 43  47
 ----- ------------------------------------------ . . .
    3 |19 21 23 27 29 33 35 39 45 47 53 57 59  63 .
    5 |   23 25 29 31 35 37 41 47 49 55 59 61  65 .
    7 |      27 31 33 37 39 43 49 51 57 61 63  67 .
   11 |         35 37 41 43 47 53 55 61 65 67  71 .
   13 |            39 43 45 49 55 57 63 67 69  73 .
   17 |               47 49 53 59 61 67 71 73  77 .
   19 |                  51 55 61 63 69 73 75  79 .
   23 |                     59 65 67 73 77 79  83 .
   29 |                        71 73 79 83 85  89 .
   31 |                           75 81 85 87  91 .
   37 |                              87 91 93  97 .
   41 |                                 95 97 101 .
   43 |                                    99 103 .
   47 |                                       107 .
      ............................................
      .                                             .
      .                                               .
      .                                                 .

 -----
|17   |
|  +  |
|     | 3  5  7 11 13 17 19 23 29 31 37 41  43  47
 ----- ------------------------------------------ . . .
    3 |23 25 27 31 33 37 39 43 49 51 57 61  63  67 .
    5 |   27 29 33 35 39 41 45 51 53 59 63  65  69 .
    7 |      31 35 37 41 43 47 53 55 61 65  67  71 .
   11 |         39 41 45 47 51 57 59 65 69  71  75 .
   13 |            43 47 49 53 59 61 67 71  73  77 .
   17 |               51 53 57 63 65 71 75  77  81 .
   19 |                  55 59 65 67 73 77  79  83 .
   23 |                     63 69 71 77 81  83  87 .
   29 |                        75 77 83 87  89  93 .
   31 |                           79 85 89  91  95 .
   37 |                              91 95  97 101 .
   41 |                                 99 101 105 .
   43 |                                    103 107 .
   47 |                                        111 .
      ............................................
      .                                             .
      .                                               .
      .                                                 .
```

```
 -----
|19   |
|  +  |
|     | 3  5  7 11 13 17 19 23 29 31 37  41  43  47
 ----- --------------------------------------------- . . .
    3 |25 27 29 33 35 39 41 45 51 53 59  63  65  69 .
    5 |   29 31 35 37 41 43 47 53 55 61  65  67  71 .
    7 |      33 37 39 43 45 49 55 57 63  67  69  73 .
   11 |         41 43 47 49 53 59 61 67  71  73  77 .
   13 |            45 49 51 55 61 63 69  73  75  79 .
   17 |               53 55 59 65 67 73  77  79  83 .
   19 |                  57 61 67 69 75  79  81  85 .
   23 |                     65 71 73 79  83  85  89 .
   29 |                        77 79 85  89  91  95 .
   31 |                           81 87  91  93  97 .
   37 |                              93  97  99 103 .
   41 |                                 101 103 107 .
   43 |                                     105 109 .
   47 |                                         113 .
      ...............................................
      .                                               .
      .                                                 .
      .                                                   .

 -----
|23   |
|  +  |
|     | 3  5  7 11 13 17 19 23 29 31 37  41  43  47
 ----- --------------------------------------------- . . .
    3 |29 31 33 37 39 43 45 49 55 57 63  67  69  73 .
    5 |   33 35 39 41 45 47 51 57 59 65  69  71  75 .
    7 |      37 41 43 47 49 53 59 61 67  71  73  77 .
   11 |         45 47 51 53 57 63 65 71  75  77  81 .
   13 |            49 53 55 59 65 67 73  77  79  83 .
   17 |               57 59 63 69 71 77  81  83  87 .
   19 |                  61 65 71 73 79  83  85  89 .
   23 |                     69 75 77 83  87  89  93 .
   29 |                        81 83 89  93  95  99 .
   31 |                           85 91  95  97 101 .
   37 |                              97 101 103 107 .
   41 |                                 105 107 111 .
   43 |                                     109 113 .
   47 |                                         117 .
      ...............................................
      .                                               .
      .                                                 .
      .                                                   .
```

```
 -----
|29    |
|  +   |
|      | 3  5  7 11 13 17 19 23 29 31  37  41  43  47
 ----- ------------------------------------------------ . . .
    3 |35 37 39 43 45 49 51 55 61 63  69  73  75  79 .
    5 |   39 41 45 47 51 53 57 63 65  71  75  77  81 .
    7 |      43 47 49 53 55 59 65 67  73  77  79  83 .
   11 |         51 53 57 59 63 69 71  77  81  83  87 .
   13 |            55 59 61 65 71 73  79  83  85  89 .
   17 |               63 65 69 75 77  83  87  89  93 .
   19 |                  67 71 77 79  85  89  91  95 .
   23 |                     75 81 83  89  93  95  99 .
   29 |                        87 89  95  99 101 105 .
   31 |                           91  97 101 103 107 .
   37 |                              103 107 109 113 .
   41 |                                  111 113 117 .
   43 |                                      115 119 .
   47 |                                          123 .
      ................................................
      .                                                 .
      .                                                   .
      .                                                     .


 -----
|31    |
|  +   |
|      | 3  5  7 11 13 17 19 23 29 31  37  41  43  47
 ----- ------------------------------------------------ . . .
    3 |37 39 41 45 47 51 53 57 63 65  71  75  77  81 .
    5 |   41 43 47 49 53 55 59 65 67  73  77  79  83 .
    7 |      45 49 51 55 57 61 67 69  75  79  81  85 .
   11 |         53 55 59 61 65 71 73  79  83  85  89 .
   13 |            57 61 63 67 73 75  81  85  87  91 .
   17 |               65 67 71 77 79  85  89  91  95 .
   19 |                  69 73 79 81  87  91  93  97 .
   23 |                     77 83 85  91  95  97 101 .
   29 |                        89 91  97 101 103 107 .
   31 |                           93  99 103 105 109 .
   37 |                              105 109 111 115 .
   41 |                                  113 115 119 .
   43 |                                      117 121 .
   47 |                                          125 .
      ................................................
      .                                                 .
      .                                                   .
      .                                                     .
```

```
 -----
|37   |
|  +  |
|     | 3  5  7 11 13 17 19 23 29 31  37  41  43  47
 ----- ------------------------------------------------ . . .
    3 |43 45 47 51 53 57 59 63 69 71  77  81  83  87 .
    5 |   47 49 53 55 59 61 65 71 73  79  83  85  89 .
    7 |      51 55 57 61 63 67 73 75  81  85  87  91 .
   11 |         59 61 65 67 71 77 79  85  89  91  95 .
   13 |            63 67 69 73 79 81  87  91  93  97 .
   17 |               71 73 77 83 85  91  95  97 101 .
   19 |                  75 79 85 87  93  97  99 103 .
   23 |                     83 89 91  97 101 103 107 .
   29 |                        95 97 103 107 109 113 .
   31 |                           99 105 109 111 115 .
   37 |                              111 115 117 121 .
   41 |                                  119 121 125 .
   43 |                                      123 127 .
   47 |                                          131 .
      ................................................
      .                                                 .
      .                                                   .
      .                                                     .


 -----
|41   |
|  +  |
|     | 3  5  7 11 13 17 19 23 29  31  37  41  43  47
 ----- ------------------------------------------------- . . .
    3 |47 49 51 55 57 61 63 67 73  75  81  85  87  91 .
    5 |   51 53 57 59 63 65 69 75  77  83  87  89  93 .
    7 |      55 59 61 65 67 71 77  79  85  89  91  95 .
   11 |         63 65 69 71 75 81  83  89  93  95  99 .
   13 |            67 71 73 77 83  85  91  95  97 101 .
   17 |               75 77 81 87  89  95  99 101 105 .
   19 |                  79 83 89  91  97 101 103 107 .
   23 |                     87 93  95 101 105 107 111 .
   29 |                        99 101 107 111 113 117 .
   31 |                           103 109 113 115 119 .
   37 |                               115 119 121 125 .
   41 |                                   123 125 129 .
   43 |                                       127 131 .
   47 |                                           135 .
      .................................................
      .                                                  .
      .                                                    .
      .                                                      .
```

| 43 + | 3 | 5 | 7 | 11 | 13 | 17 | 19 | 23 | 29 | 31 | 37 | 41 | 43 | 47 | . . . |
|---|---|---|---|---|---|---|---|---|---|---|---|---|---|---|---|
| 3 | 49 | 51 | 53 | 57 | 59 | 63 | 65 | 69 | 75 | 77 | 83 | 87 | 89 | 93 | . |
| 5 | | 53 | 55 | 59 | 61 | 65 | 67 | 71 | 77 | 79 | 85 | 89 | 91 | 95 | . |
| 7 | | | 57 | 61 | 63 | 67 | 69 | 73 | 79 | 81 | 87 | 91 | 93 | 97 | . |
| 11 | | | | 65 | 67 | 71 | 73 | 77 | 83 | 85 | 91 | 95 | 97 | 101 | . |
| 13 | | | | | 69 | 73 | 75 | 79 | 85 | 87 | 93 | 97 | 99 | 103 | . |
| 17 | | | | | | 77 | 79 | 83 | 89 | 91 | 97 | 101 | 103 | 107 | . |
| 19 | | | | | | | 81 | 85 | 91 | 93 | 99 | 103 | 105 | 109 | . |
| 23 | | | | | | | | 89 | 95 | 97 | 103 | 107 | 109 | 113 | . |
| 29 | | | | | | | | | 101 | 103 | 109 | 113 | 115 | 119 | . |
| 31 | | | | | | | | | | 105 | 111 | 115 | 117 | 121 | . |
| 37 | | | | | | | | | | | 117 | 121 | 123 | 127 | . |
| 41 | | | | | | | | | | | | 125 | 127 | 131 | . |
| 43 | | | | | | | | | | | | | 129 | 133 | . |
| 47 | | | | | | | | | | | | | | 137 | . |

| 47 + | 3 | 5 | 7 | 11 | 13 | 17 | 19 | 23 | 29 | 31 | 37 | 41 | 43 | 47 | . . . |
|---|---|---|---|---|---|---|---|---|---|---|---|---|---|---|---|
| 3 | 53 | 55 | 57 | 61 | 63 | 67 | 69 | 73 | 79 | 81 | 87 | 91 | 93 | 97 | . |
| 5 | | 57 | 59 | 63 | 65 | 69 | 71 | 75 | 81 | 83 | 89 | 93 | 95 | 99 | . |
| 7 | | | 61 | 65 | 67 | 71 | 73 | 77 | 83 | 85 | 91 | 95 | 97 | 101 | . |
| 11 | | | | 69 | 71 | 75 | 77 | 81 | 87 | 89 | 95 | 99 | 101 | 105 | . |
| 13 | | | | | 73 | 77 | 79 | 83 | 89 | 91 | 97 | 101 | 103 | 107 | . |
| 17 | | | | | | 81 | 83 | 87 | 93 | 95 | 101 | 105 | 107 | 111 | . |
| 19 | | | | | | | 85 | 89 | 95 | 97 | 103 | 107 | 109 | 113 | . |
| 23 | | | | | | | | 93 | 99 | 101 | 107 | 111 | 113 | 117 | . |
| 29 | | | | | | | | | 105 | 107 | 113 | 117 | 119 | 123 | . |
| 31 | | | | | | | | | | 109 | 115 | 119 | 121 | 125 | . |
| 37 | | | | | | | | | | | 121 | 125 | 127 | 131 | . |
| 41 | | | | | | | | | | | | 129 | 131 | 135 | . |
| 43 | | | | | | | | | | | | | 133 | 137 | . |
| 47 | | | | | | | | | | | | | | 141 | . |

99)Smarandache-Vinogradov sequence:

0,0,0,0,1,2,4,4,6,7,9,10,11,15,17,16,19,19,23,25,26,26,28,33,32,35,43,39,
40,43,43,...

(a(2k+1) represents the number of different combinations such that 2k+1 is written as a sum of three odd primes.)

This sequence is deduced from the Smarandache-Vinogradov table.

100)Circular sequence:

```
1,12,21,123,231,312,1234,2341,3412,4123,12345,23451,34512,45123,51234,
| |   | |         | |                 | |                           |
   ---   ---------   -----------------   ---------------------------
1   2        3               4                         5


123456,234561,345612,456123,561234,612345,1234567,2345671,3456712,...
|                                       | |
 ---------------------------------------   ----------------------- ...
                    6                                          7
```

101)Simple numbers:

2,3,4,5,6,7,8,9,10,11,13,14,15,17,19,21,22,23,25,26,27,29,31,33,34,35,37,38,
39,41,43,45,46,47,49,51,53,55,57,58,61,62,65,67,69,71,73,74,77,78,79,82,83,
85,86,87,89,91,93,94,95,97,101,103,...

(A number n is called simple number if the product of its proper divisors is less than or equal to n.)

Generally speaking, n has the form:

n = p, or p^2, or p^3, or pq,  where p and q are distinct primes.